\documentclass[12pt]{amsart}
\usepackage{amssymb}

\usepackage[all]{xy}
\newtheorem{theorem}{Theorem}[section]

\newtheorem{proposition}[theorem]{Proposition}
\newtheorem{corollary}[theorem]{Corollary}

\theoremstyle{definition}
\newtheorem{definition}[theorem]{Definition}
\newtheorem{example}[theorem]{Example}

\theoremstyle{remark}
\newtheorem{remark}[theorem]{Remark}

\numberwithin{equation}{section}

\begin{document}

\title[Totally bounded sets in the absolute weak topology]{Totally bounded sets in the absolute weak topology}

\author[H. Ardakani]
{Halimeh Ardakani }

\address{Department of Mathematics, Payame Noor University, Tehran, Iran}
 \email{halimeh\_ardakani@yahoo.com, ardakani@pnu.ac.ir}

\author[J. X. Chen]
{Jin Xi Chen}
\address{School of Mathematics, Southwest Jiaotong University, Chengdu 610031, China}
\email{jinxichen@swjtu.edu.cn}

%    General info
\subjclass[2010]{Primary 46B42; Secondary 46B50, 47B65}
%\dedicatory{}

\keywords{almost $L$-set, $|\sigma|(E,E^\prime)$-totally bounded set, PL-compact operator, almost Dunford-Pettis operator, Banach lattice.}

\begin{abstract}
		
	In this paper, almost Dunford-Pettis operators with ranges in $c_0$ are used to identify totally bounded sets in the absolute weak topology. That is, a bounded subset $A$ of a Banach lattice $E$ is $|\sigma|(E,E^\prime)$-totally bounded if and only if  $T(A)\subset c_0$ is relatively compact for every almost Dunford-Pettis operator $T:E\to{c_0}$. As an application, we show that for two Banach lattices $E$ and $F$ every positive operator from $E$ to $F$ dominated by a PL-compact operator is PL-compact if and only if either the norm of $E^{\,\prime}$ is order continuous  or  every order interval in $F$ is $|\sigma|(F,F^\prime)$-totally bounded.

\end{abstract}

\maketitle\baselineskip 4.9mm
	
	\section{Introduction and preliminaries}

 Throughout this paper, we denote Banach spaces  by $X, Y$, and denote Banach lattices by $E, F$. $B_X$ is the closed unit ball of $X$. $E^{+}$ denotes the positive cone of $E$ and $Sol(A)$ denotes the solid hull of a subset $A$ of $E$.

\par A bounded subset $A$ of  $X$ is called a \textit{limited } (resp. \textit{Dunford-Pettis}) \textit{set} if $\sup_{x\in{A}}{\vert{f_{n}(x)}\vert}\to{0}$ for  each sequence $(f_{n})$ of $ X^{\prime}$ satisfying $f_{n}\xrightarrow{w^*}0$ (resp. $f_{n}\xrightarrow{w}0$),  or equivalently, if  every bounded (resp. weakly compact) operator from $X$ to $c_0$ carries $A$ to a relatively compact set.  See, e.g., \cite{Positive,Andrews,Bourgain}. Dually, a bounded subset $B$ of $X^{\prime}$ is called an \textit{$L$\,-set} if $\sup_{f\in{B}}{\vert{f(x_n)}\vert}\to{0}$ for each sequence $(x_n)$ in $X$ with $x_{n}\xrightarrow{w}0$ (see, e.g., \cite{L-set,Bator,Em}). In particular, we say a sequence $(f_n)$ in $X^{\prime}$ is an $L$\,-sequence if $\{f_n:n\in \mathbb{N} \}$ is an $L$\,-set. It is easily verified that, to each Dunford-Pettis operator $T:X\to{c_0}$, there corresponds a unique weak$^{*}$-null $L$\,-sequence $(f_n)$ in $X^{\prime}$ such that $Tx=(f_{n}(x))_{n=1}^{\infty}$.

\par Let us recall that a  subset $A$ of  $X$ is said to be \textit{weakly  precompact} or \textit{conditionally weakly compact} if every sequence in $A$ has a weak Cauchy subsequence.  Based on the work of Odell and Stegall, Rosenthal \cite[p.377]{1977Point} gave an operator characterization of weakly precompact sets: \textit{a subset $A$ of  $X$ is weakly precompact if and only if for every Banach space $Y$ and
every Dunford-Pettis operator $T:X\to{Y}$, $T(A)$ is  relatively compact}.  In her paper \cite{Ioana}, Ghenciu  proved that $A\subset X$ is  weakly precompact if and only if $T(A)$ is relatively compact for every Dunford-Pettis operator $T:X\to{c_0}$ \cite[Theorem 1\&Corollary 9]{Ioana}. Thus, Dunford-Pettis operators with ranges in $c_0$ can also be employed to identify weakly precompact sets. Based upon this characterization of weakly precompact sets,  Xiang, Chen and Li \cite{XCL} recently considered weak precompactness properties in Banach lattices.

\par To justify our motivation to present this work, we have to mention a class of operators which is much like Dunford-Pettis operators. Let us recall that a bounded linear operator $T: E\rightarrow X$ from a Banach lattice to a Banach space is called an \textit{almost Dunford-Pettis operator} if $\|Tx_{n}\| \rightarrow 0$ for each disjoint weakly null sequence $(x_n)$ in $E$, or equivalently, if $\|Tx_{n}\| \rightarrow 0$ whenever $0\leq x_{n}\xrightarrow{w}0$ in $E$ \cite{some,Sanchez}. From \cite[Example 4), p.230]{Wnuk} it follows that every positive almost Dunford-Pettis operator from $E$ to $c_0$ is Dunford-Pettis. In general, an almost Dunford-Pettis operator is not necessarily Dunford-Pettis.  For instance, every bounded linear operator from $L^{1}[0,1]$ to $c_0$ is almost Dunford-Pettis since $L^{1}[0,1]$ has the positive Schur property. However,  the  operator $T:L^{1}[0,1]\to{c_0}$ defined by $$Tf=\left(\int_{0}^{1}f(t)r_{1}(t)\textrm{d}t, \int_{0}^{1}f(t)r_{2}(t)\textrm{d}t, \cdot\cdot\cdot, \int_{0}^{1}f(t)r_{n}(t)\textrm{d}t, \cdot\cdot\cdot\right),$$where  $(r_n)$ is the sequence of Rademacher functions, is not Dunford-Pettis since $r_{n}\xrightarrow{w}{0}$ in $L^{1}[0,1]$ and $\|Tr_n\|=1$ for all $n$. 

\par Let $A$ be a bounded subset of a Banach lattice $E$. One may ask the following question:\\

\textit{What properties does the set $A$ possess if $T(A)$ is relatively compact in $c_0$ for each almost Dunford-Pettis operator $T:E\to c_0$?}\\

\par The  first objective of this paper is to introduce and study a class of sets that we  call super weakly precompact sets for the time being, which can be identified by almost Dunford-Pettis operators with range in $c_0$. That is, a bounded subset $A$ of $E$  is called \textit{super weakly precompact} if $T(A)\subset c_0$ is relatively compact for every almost Dunford-Pettis operator $T:E\to{c_0}$. Clearly, limited sets, disjoint weakly null sequences and  positive weakly null sequences  in $E$ are all super weakly precompact. Also, every super weakly precompact set is weakly precompact. On the other hand, as we have mentioned above,  the sequence $(r_n)$ of Rademacher functions is a relatively weakly compact subset of $L^{1}[0,1]$. However, it is not super weakly precompact. By definition, the closed absolutely convex hull of a super weakly precompact set is likewise super weakly precompact.
 \par With the progress of our study we found that every positive operator from $E$ to an $L$-space carries super weakly precompact subsets of $E$ onto relatively compact sets. This property was originally enjoyed by the class of $|\sigma|(E,E^\prime)$-totally bounded sets in $E$.  Inspired by this observation, we conjecture that the so-called super weakly precompact sets can be closely related to totally bounded  sets with respect the absolute weak topology. Recall that the absolute weak topology $|\sigma|(E,E^\prime)$ on $E$ is a locally convex-solid topology generated by the family of lattice seminorms $\{p_f:f\in E^\prime\}$, where $p_f$ is defined via the formula $p_{f}(x):=|f|(|x|)$ for $x\in E$. Then, a subset $A$ of $E$ is $|\sigma|(E,E^\prime)$--totally bounded if for every $\varepsilon>0$ and every finite collection $\{f_{1}, f_{2},\cdots,f_{n}\}\subset E^\prime$ there exists a finite subset $\Phi$ of $A$ such that $A\subset \Phi+\bigcap_{i=1}^{n}\{x:|f_{i}|(|x|)<\varepsilon\}$. Our second objective  is to show that a subset of $E$ is super weakly precompact if and if it is  $|\sigma|(E,E^\prime)$--totally bounded (Theorem \ref{wo de PL-compact}). Hence, we obtain an alternative characterization of $|\sigma|(E,E^\prime)$--totally bounded sets: $A\subset E$ is $|\sigma|(E,E^\prime)$--totally bounded whenever $T(A)\subset c_0$ is relatively compact for every almost Dunford-Pettis operator $T:E\to{c_0}$.
 
\par  Finally, we consider the operators related to totally bounded sets in absolute weak topologies. Following Dodds and Fremlin \cite[Definition 4.1]{Dodds and Fremlin PL-compact},  an operator $T:X\to{E}$ is called a \textit{PL-compact operator} whenever $TB_X$ is $|\sigma|(E,E^\prime)$--totally bounded in $E$.  Theorem 4.6 of \cite{Dodds and Fremlin PL-compact} asserts that, for two Banach lattices $E$ and $F$ such that $E^{\,\prime}$ and $F$ have order continuous norms, the set of regular PL-compact operators from $E$ to $F$ forms a band in $\mathcal{L}^{r}(E,F)$.   As an application of our results, we  solve the domination problem for positive PL-compact operators. we show that for Banach lattices $E$ and $F$ every positive operator from $E$ to $F$ dominated by a PL-compact operator is PL-compact if and only if either the norm of $E^{\,\prime}$ is order continuous  or  every order interval in $F$ is $|\sigma|(F,F^\prime)$--totally bounded ( Theorem \ref{domination}).

\par For our further discussion, we shall need some notions from Banach lattice theory.
\begin{enumerate}
  \item A bounded subset $A$ of a Banach lattice $E$ is called
 \begin{itemize}
   \item[$\raisebox{0mm}{------}$]   an  \textit{almost Dunford-Pettis set}  if every weakly null disjoint sequence $(f_{n})$ of $ E^{\,\prime}$ converges uniformly to zero on $A$ \cite{Almost DP set}.
   \item[$\raisebox{0mm}{------}$] a \textit{positively limited set} if every positive weak$^*$-null
sequence $(f_n)$ in $(E^{\,\prime})^+$ converges uniformly to zero on $A$ \cite{ positively limited set}.
\item[$\raisebox{0mm}{------}$] an \textit{$L$-weakly compact set} if  $\|x_{n}\|\to 0$ for each disjoint sequence $(x_n)$ in $Sol(A)$.
       \end{itemize}
  \item A bounded subset $B\subset E^{\,\prime}$ is called an \textit{almost L-set} if every disjoint weakly null sequence $(x_n)$ in $E$  converges to zero uniformly on $B$ \cite{L-set}. A sequence $(f_{n})$ of $ E^{\,\prime}$ will be called an almost $L$-sequence whenever $\{f_{n}:n\in \mathbb{N}\}$ is an almost $L$-set.

  \item A Banach lattice $E$ has the \textit{weak Dunford-Pettis property} (resp. \textit{positive DP$^{\,*}$ property}) if every relatively weakly compact set in $E$ is almost Dunford-Pettis (resp. positively limited).
   \item A Banach lattice $E$ has the \textit{positive Schur property}  if every disjoint weakly null sequence in $E$ is norm null, or equivalently, if $\|x_n\|\to 0$ whenever $0\leq x_{n}\xrightarrow{w}0$ in $E$. Every $L$-space has the positive Schur property, and every Banach lattice with the positive Schur property is a KB-space.
\end{enumerate}
\par The definitions and notions from Banach lattice theory which appear here, are standard. We refer the reader to the references \cite{Locally Solid,Positive,Meyer}. 

\section{Super weakly precompact sets and $|\sigma|(E,E^\prime)$--totally bounded sets}

The super weakly precompact sets can  be identified by employing  almost $L$-sequences which are weak$^*$-null.

\begin{proposition}\label{d} For a bounded subset $A$ of a Banach lattice $E$ the following assertions are equivalent:
	\begin{itemize}
		\item [(a)] $A$ is a super weakly precompact.
        \item [(b)] Every  weak$^{\,*}$-null almost $L$-sequence
	 $(f_n)$ in $E^\prime$ converges uniformly to zero on $A$, that is, $\sup_{x\in A}|f_{n}(x)|\rightarrow0$ \,\,$(n\rightarrow\infty)$.
		\item [(c)] $f_{n}(x_n)\rightarrow0$  for every weak$^{\,*}$-null almost $L$-sequence
		 $(f_n)$ in $E^\prime$ and every sequence $(x_n)\subset A$.
	\end{itemize}
\end{proposition}
\begin{proof}
$(a)\Leftrightarrow(b).$ Since each almost Dunford-Pettis operator $T:E \rightarrow c_0$ is uniquely determined by a weak$^*$-null almost $L$-sequence $(f_n)$ in $E^\prime$ such that $T(x)=(f_n(x))$ for all $x \in E$ and also, the subset $T(A)$ of $c_0$ is relatively compact if and only if $$s_n=\sup_{(b_{n})_{n=1}^{\infty}\in T(A)}|b_n|= \sup_{x\in A}| f_n (x)| \rightarrow 0,$$ the desired result can be easily proved.\\
 $(b)\Rightarrow (c).$ Obvious.\\
	$(c)\Rightarrow (b).$ Let $(f_n)$ be a weak$^*$-null almost $L$-sequence in $E^\prime$. For each $n$, we can choose $x_n$ in $A$ such that $\sup_{x\in A}|f_n(x)|<|f_n(x_n)|+\frac{1}{n}$. By hypothesis (c), it follows that $\sup_{x\in A}|f_{n}(x)|\rightarrow0$ \,\,$(n\rightarrow\infty)$.
\end{proof}

\par As a practical result throughout this paper, the following result tells us when every bounded subset of a Banach lattice is super weakly precompact.

 \begin{theorem}\label{k.115} For a Banach lattice $E$ the following assertions are equivalent:
 	\begin{itemize}
 		\item [(a)] $B_E$ is a super weakly precompact set, that is, every almost Dunford-Pettis operator $T:E\to{c_0}$ is compact.

 		\item [(b)] $E^\prime$ is a discrete Banach lattice with order continuous norm.
       \item [(c)] Every almost Dunford-Pettis operator from $E$ to an arbitrary Banach space is compact.
 	\end{itemize}
 \end{theorem}
 \begin{proof}
 	$(a)\Rightarrow (b).$ It suffices to show that $[-f,f]$ is compact for each $0\leq f\in E^\prime$.
 	Since $B_E$ is super weakly precompact (and hence, $B_E$ is weakly precompact), $E$ contains no isomorphic copy of $\ell_1$ and so, $E^\prime$ has order continuous norm (cf. \cite[Theorem 4.56 \& 4.69]{Positive}).  Hence, $[-f,f]$ is weakly compact. Assume by way of contradiction that $[-f,f]$ is not compact. Then there would exist a weakly convergent sequence $(f_n)\subset[-f,f]$ which contains no norm convergent subsequences. Let $f_0$ be the weak limit of $(f_n)_{n=1}^{\infty}$. Then $(f_{n}-f_0)_{n=1}^{\infty}$ is weakly null and order bounded. Note that each order interval in a dual of a Banach lattice is an almost $L$-set. Therefore,  $(f_{n}-f_0)_{n=1}^{\infty}$ is also an almost $L$-sequence. From the super weak precompactness of $B_E$  and Proposition \ref{d} it follows that $\|f_{n}-f_0\|\to{0}$. This leads to a contradiction.
 \par $(b)\Rightarrow (a).$ Let $(f_n)$ be a weak$^*$-null almost $L$-sequence in $E^\prime$. Since $E^\prime$ has order continuous norm, the sequence $(f_n)$ is $L$-weakly compact \cite[Lemma 2.9]{Chen}. Since  $E^\prime$ is discrete and has order continuous norm,  it follows that $(f_n)$ is relatively compact. Hence $\|f_{n}\|\rightarrow 0$; that is, $B_E$ is super weakly precompact.
 \par $(a), (b)\Rightarrow (c).$ Let $T:E\to{X}$ be an almost Dunford-Pettis operator from $E$ to a Banach space $X$. We first prove that $T$ is a Dunford-Pettis operator. To this end, let $x_{n}\xrightarrow{w}0$ in $E$ and $0\leq f\in E^{\prime}$. Then, $f(|x_n|)=\sup_{g\in [-f,f]}|g(x_n)|\to 0$ since $E^\prime$ is  discrete  with order continuous norm and hence $[-f,f]$ is compact. This implies that $|x_n|\xrightarrow{w}0$ (and hence, both $(x_{n}^+)$ and $(x_{n}^-)$ are weakly null). From the almost Dunford-Pettis property of $T$ it follows that $\|Tx_n\|\leq\|T(x_{n}^+)\|+\|T(x_{n}^-)\|\to0$. That is, $T$ is Dunford-Pettis. Since $B_E$ is super weakly precompact and hence $B_E$ is  weakly precompact, $T(B_E)$ is relatively compact.
 \par $(c)\Rightarrow (a).$ Obvious.
 \end{proof}

\par Recall that  $E$ has the weak Dunford-Pettis property  if and only if every weakly compact operator from $E$ to $c_0$ is almost Dunford-Pettis \cite{Wnuk}. The following characterization of the Schur property improves \cite[Theorem 2.10]{Schur}. Here we give a direct and simpler proof.
 \begin{corollary}\label{tlh2.1} For a Banach lattice $E$ the following statements are equivalent:
 	\begin{itemize}
 		\item [(a)] $E^\prime$ has the Schur property.
 		\item [(b)] $E$ has the weak Dunford-Pettis property and $E^\prime$ is discrete with order continuous norm.
 	\end{itemize}	
 \end{corollary}
 \begin{proof}
 	$(a) \Rightarrow (b)$. The Schur property of $E^{\prime}$ implies that $E^\prime$ has the Dunford-Pettis property  and its order intervals are all compact. Therefore, $E$ has the weak Dunford-Pettis property and $E^\prime$ is discrete with order continuous norm.
	\par $(b) \Rightarrow (a)$. Let $(f_n)\subset E^{\prime}$ satisfy $f_{n}\xrightarrow{w}0$ in $E^{\prime}$. Then the operator $T:E\to{c_0}$, where $Tx=(f_{1}(x), f_{2}(x), \ldots)$, is weakly compact. Since $E$ has the weak Dunford-Pettis property, $T$ is  almost Dunford-Pettis. On the other hand, since $E^{\prime}$ is discrete with order continuous norm, fom Theorem \ref{k.115} it follows that $T(B_E)$ is a relatively compact subset of $c_0$.  Thus, $\|f_n\|=\sup_{x\in B_E}|f_{n}(x)|\to{0}$.  This implies that $E^{\prime}$ has the Schur property.
 \end{proof}

 \par It is clear that each limited set in a Banach lattice is super weakly precompact. Conversely, a super weakly precompact set is not necessarily limited. By Theorem \ref{k.115}, $B_{c_0}$ is indeed  super weakly precompact  in $c_0$. However, by the well-known Josefson--Nissenzweig theorem, $B_X$ can not be a limited subset of $X$ whenever the Banach space $X$ is infinite dimensional (see, e.g., \cite[p. 219]{Diestel}). We know that  disjoint weakly null sequences are all super weakly precompact. The following result shows that the limitedness of super weakly precompact sets is determined by the limitedness of disjoint weakly null sequences.

 \begin{theorem}\label{2gy.3} For a Banach lattice $E$ the following statements are equivalent:
 	\begin{itemize}
 		\item [(a)] Every super weakly precompact subset of $E$ is a limited set.
 		\item [(b)] Every disjoint weakly null sequence in $E$ is limited.
 		\item [(c)] Every weak$^*$-null sequence in $E^{\prime}$ is an almost $L$-sequence.
 	\end{itemize}
 \end{theorem}

 \begin{proof}
 The implications of $(a)\Rightarrow (b)$ and  $(c)\Rightarrow (a)$ are obvious.
 \par $(b\Rightarrow (c)$ Let $(f_n)$ be a weak$^*$-null sequence in $E^{\prime}$. For  every disjoint weakly null sequence $(x_n) \subset E$, by our hypothesis, $(x_n)$ is limited.  Therefore, $f_n(x_n)\rightarrow 0$. This implies that $(f_n)$ be a an almost $L$-sequence in $E^{\prime}$.
  \end{proof}
\par Recall that a Banach space $X$ called a \textit{Gelfand--Phillips space} if every  limited set in $X$ is relatively compact.  It  should be noted that a $\sigma$-Dedekind complete Banach lattice  is a Gelfand--Phillips space if and only if its norm is order continuous. See, e.g., \cite[Theorem 4.5, p.80]{WOrder}. Also, a Banach lattice $E$ has the positive Schur property if and only if every relatively weakly compact subset of $E$ is $L$-weakly compact.
\begin{corollary}\label{u2y.3} For a Banach lattice $E$ the following assertions are equivalent:
\begin{itemize}
		\item [(a)] $E$ has the positive Schur property.
        \item [(b)]	Each super weakly precompact set in $E$ is relatively compact.
        \item [(c)] Each super weakly precompact set in $E$ is $L$-weakly compact.
\end{itemize}
\end{corollary}

\begin{proof}
$(a)\Rightarrow (b)$. The positive Schur property of $E$ implies that every disjoint weakly null sequence in $E$ is norm null.  By Theorem \ref{2gy.3}, every super weakly precompact subset of $E$ is a limited set, and hence it is relatively compact since $E$ is a Gelfand--Phillips space.\\
$(b)\Rightarrow (a)$. It follows easily from the fact that every disjoint weakly null sequence in $E$ is super weakly precompact.\\
$(a),(b)\Rightarrow (c)$. Let $A$ be a super weakly precompact subset of $E$. Then by our hypothesis, $A$ is relatively compact. Therefore, $A$ is $L$-weakly compact since $E$ has order continuous norm.\\
$(c)\Rightarrow (a)$. Let $(x_n)\subset E$ be a disjoint sequence such that $x_{n}\xrightarrow{w}0$. By our hypothesis, $(x_n)$ is $L$-weakly compact. Therefore, $\|x_n\|\to0$. That is, $E$ has the positive Schur property.
\end{proof}

\par As we have pointed out in the introduction, super weak precompactness implies weak precompactness whereas even a weakly compact set is not necessarily  super weak precompact. The following result tells us under what conditions weak precompactness and super weak precompactness  coincide.

\begin{theorem}\label{2y7.1} For a Banach lattice $E$ the following assertions are equivalent:
	\begin{itemize}
		\item [(a)] Every weakly precompact subset of $E$ is super weakly precompact.
		\item [(b)] Every relatively weakly compact subset of $E$ is super weakly precompact.
		\item [(c)] $E$ has  weakly sequentially continuous lattice operations (i.e., $x_{n}\xrightarrow{w}{0}$ in $E$ implies $|x_n|\xrightarrow{w}0$ in $E$).
        \item [(d)] Every almost Dunford-Pettis operator from $E$ to an arbitrary Banach space is Dunford-Pettis.
	\end{itemize}
\end{theorem}
\begin{proof}
	It suffices to prove $(b)\Rightarrow (c).$ We assume by way of contradiction that the lattice operations in $E$ are not weakly sequentially continuous. Then there would exist a weakly null sequence $(x_n)\subset E$ such that $(|x_{n}|)$ is not weakly null. Hence there exists  $0\leq f\in E^\prime$ satisfying  $f(|x_n|)>\varepsilon$ for some $\varepsilon>0$ and for all $n$. As in the proof of Theorem 2 of \cite{converse}, we can find a sequence $(h_n)\subset [-f,f]$ such that $h_{n}\xrightarrow{w^*}0$ and $h_n(x_n)\ge \varepsilon$ for all $n$. Note that $(h_n)$, is also an almost $L$-sequence since every order bounded set in the dual of a Banach lattice is  an almost $L$-set. In view of Theorem \ref{d}, this is impossible.
\end{proof}
\par We know that bounded linear operators can preserve many  topological properties. However, the following example shows that they do not necessarily do the same thing to super weakly precompact sets.
\begin{example}\label{t.1} Let $T:\ell_{2} \to L^{1}[0,1]$ be the isomorphic embedding of $\ell_{2}$ in $L^{1}[0,1]$.  Khinchine's inequality implies that such an  isomorphism certainly exists (see, e.g., \cite[Theorem 2.25]{Ryan}).
By Theorem \ref{k.115},  $B_{\ell_{2}}$ is a super weakly precompact subset of $\ell_2$.  Assume to the contrary that $TB_{\ell_{2}}$ is super weakly precompact in $L^1[0,1]$. Then, by Corollary \ref{u2y.3}, $TB_{\ell_{2}}$ is a compact subset of $L^1[0,1]$. This is impossible.
\end{example}

\par The order bounded operators can indeed preserve the super weak precompactness property.

 \begin{proposition}\label{image under order bounded operator} Let $T:E\to{F}$ be an order bounded operator between Banach lattices.
 \begin{enumerate}
   \item If $A$ is a super weakly precompact set in $E$, then $T(A)$ is likewise a super weakly precompact set in $F$.
   \item If $B$ is an almost $L$-set in $F^{\,\prime}$, then $T^{\,\prime}(B)$ is likewise an almost $L$-set in $E^{\,\prime}$.
 \end{enumerate}
\end{proposition}
\begin{proof}
(1) It suffices to prove that  $S(TA)$ is a relatively compact subset of $c_0$ for  an arbitrary almost Dunford-Pettis oprator $S:F\to{c_0}$. First we claim that $ST:E\to{c_0}$ is also almost Dunford-Pettis. To this end, let $(x_n)$ be a disjoint sequence in $E$ such that $ x_{n}\xrightarrow{w}0$. Then, for  $0\leq f\in F^{\,\prime}$, by the Riesz-Kantorovich Formula we have
$$\langle f, |Tx_n|\rangle=\sup_{g\in [-f,f]}|g(Tx_n)|=\sup_{g\in [-f,f]}|\langle T^{\,\prime}g, x_n\rangle|=\sup_{h\in T^{\,\prime}[-f,f]}| h (x_n)|\xrightarrow{n\rightarrow\infty}{0}$$since $T^{\,\prime}$ is also an order bounded operator and hence $T^{\,\prime}[-f,f]$ is an almost $L$-set in $E^{\,\prime}$ \cite[Proposition 3.4]{L-set}. This implies that $|Tx_n|\xrightarrow{w}{0 }$. It follows that $\|S(Tx_n)\|\to{0}$, since $S$ is an almost Dunford-Pettis operator. That is, $ST:E\to{c_0}$ is  almost Dunford-Pettis. It follows that $ST(A)$ is relatively compact in $c_0$ since $A$ is super weakly precompact. This implies that $T(A)$ is super weakly precompact in $F$.

\par (2) Let $(x_n)\subset E$ be a disjoint sequence  such that $ x_{n}\xrightarrow{w}0$. Then, from the proof of Part (1), we can see that $|Tx_n|\xrightarrow{w}{0}$. Since $B$ be an almost $L$-subset of $F^{\,\prime}$, from \cite[Theorem 2.14]{positively limited set} it follows that
 $$\sup_{g\in B}|\langle T^{\,\prime}g, x_n\rangle|=\sup_{g\in B}|\langle g, Tx_{n}\rangle|\leq\sup_{g\in B}|\langle g, (Tx_{n})^+\rangle|+\sup_{g\in B}|\langle g, (Tx_{n})^-\rangle|\to{0}$$This implies that $T^{\,\prime}(B)$ is also an almost $L$-set in $E^{\,\prime}$.
\end{proof}

\begin{corollary}\label{2ghy.3} For a Banach lattice $E$ the following statements are equivalent:
	\begin{itemize}
		\item [(a)] Each super weakly precompact set in $E$ is relatively weakly compact.
		\item [(b)]	$E$ is weakly sequentially complete, i.e., $E$ is a KB-space.
	\end{itemize}
\end{corollary}
\begin{proof} It suffices to prove that $(a)\Rightarrow(b)$. Assume to the contrary that $E$ is not a KB-space. Then $c_0$  lattice embeds into $E$. Let $T:c_0\to E$ be the (into) lattice isomorphism. From Proposition \ref{image under order bounded operator} it follows that $TB_{c_{0}}$ is super weakly precompact in $E$ and hence, by our hypothesis $(a)$, $TB_{c_{0}}$ is relatively weakly compact and hence the operator $T$ is weakly compact. This leads to a contradiction.
\end{proof}

As we have promised in the Introduction, we now have to show that  super weak precompact sets and  $|\sigma|(E,E^\prime)$--totally bounded sets are the same.  Now it is good timing. Note that, from Corollary \ref{u2y.3} and Proposition \ref{image under order bounded operator} it follows that every positive operator from $E$ to an $L$-space maps super weakly precompact subsets of $E$ to relatively compact sets. This property was originally enjoyed by the earlier known class of PL-compact sets which was defined by P.G. Dodds and D.H. Fremlin in their excellent and classical paper \cite{Dodds and Fremlin PL-compact}. We follow the notation used in \cite{Dodds and Fremlin PL-compact}. Let $0\leq g\in E^{\,\prime}$ and $N_{g}=\{x\in E:g(|x|)=0\}$. Denote by $j_{g}$ the quotient map from $E$ onto $E/N_{g}$. $(E;g)$ is the completion of $E/N_{g}$ with respect to norm $\|[x]\|=\|j_{g}(x)\|=g(|x|)$ \,\,($x\in E$). It should be noted that $j_{g}$ is a lattice homomorphism and $(E;g)$ is an $L$-space. Following Dodds and Fremlin \cite[Definition 4.1]{Dodds and Fremlin PL-compact}, a set $A\subset E$ is called a \textit{PL-compact set} if $j_{g}(A)$ is relatively compact in $(E;g)$ for each $0\leq g\in E^{\,\prime}$.  However, now most authors commonly use the name \textit{$|\sigma|(E,E^\prime)$--totally bounded sets} in place of PL-compact sets in $E$ due to the following observations. For every finite collection $\{f_{1}, f_{2},\cdots,f_{n}\}\subset E^\prime$, we have $g=\Sigma_{i=1}^{n}|f_i|\in E^\prime$. Since $\{x:(\Sigma_{i=1}^{n}|f_{i}|)(|x|)<\varepsilon\}\subset\bigcap_{i=1}^{n}\{x:|f_{i}|(|x|)<\varepsilon\}$, we can easily see that $A$ is $|\sigma|(E,E^\prime)$--totally bounded if and only if $j_{g}(A)$ is norm totally bounded in $(E;g)$ for each $0\leq g\in E^\prime$, i.e., $A$ is PL-compact. Following \cite{WOrder}, we say a subset $A$ of $E$ is  \textit{disjointly weakly compact} whenever  $x_{n}\xrightarrow{w}0$ for every disjoint sequence   $(x_n)\subset Sol(A)$.

\begin{theorem}\label{wo de PL-compact}
A bounded subset of $E$ is super weakly precompact if and only if it is $|\sigma|(E,E^\prime)$--totally bounded.
\end{theorem}
\begin{proof}
Assume first that $A\subset E$ is super weakly precompact. For every $g\in (E^{\,\prime})^+$, by Proposition \ref{image under order bounded operator} $j_{g}(A)$ is still super weakly precompact in $(E;g)$. Since $(E;g)$ is an $L$-space (hence has the positive Schur property), from  Corollary \ref{u2y.3} it follows that $j_{g}(A)$ is a relatively compact set in $(E;g)$. This implies that $A$ is $|\sigma|(E,E^\prime)$--totally bounded in $E$.

\par For the converse, we assume that $A$ is $|\sigma|(E,E^\prime)$--totally bounded in $E$. First, we claim that $A$ is disjointly weakly compact. To this end, let $(x_n)\subset Sol(A)$ be a  disjoint sequence and let $0\leq g\in E^{\,\prime}$. Clearly, $(j_{g}(x_n))$ is a disjoint sequence in $Sol(j_{g}(A))$. Since $j_{g}(A)$ is relatively compact in the $L$-space $(E;g)$, $j_{g}(A)$ is $L$-weakly compact. Therefore, we have
    $$|g(x_n)|\leq g(|x_n|)=\|j_{g}(x_n)\|\to{0},$$That is, $x_{n}\xrightarrow{w}0$, and hence $A$ is disjointly weakly compact.
\par Now, let $(f_n)_{n=1}^{\infty}\subset F^{\prime}$ be an almost $L$-sequence satisfying $f_{n}\xrightarrow{w^*}0$. Then the solid hull $B=Sol\{f_{n}: n\in \mathbb{N}\}$ is also an almost $L$-set (cf. \cite[Theorem 3.11]{solid hull of almost L}). The disjointly weak compactness of $Sol(A)$ implies that every disjoint sequence $(x_n)$ in $Sol(A)$ is weakly null, and hence $(x_n)$ converges to zero uniformly on $B$. Therefore, for every $\varepsilon>0$ there exists $g\in (E^\prime)^+$ satisfying $$\left(|f|-g\right)^{+}(|x|)\leq\frac{\varepsilon}{3},\,\,\,\,\,\forall x\in Sol(A),\,\,\,\forall f\in B$$See, e.g., \cite[Theorem 2.3.3]{Meyer}. Since  $A$ is  $|\sigma|(E,E^\prime)$--totally bounded in $E$, there exists a finite collection $\{x_1, x_2, \ldots, x_m\}$ of elements of $A$ such that for each $x\in A$ we can find some $x_i$ satisfying $g(|x-x_i|)\leq\frac{\varepsilon}{3}.$ Consequently, from the identity $|f|=|f|\wedge g+(|f|-g)^+$ we see that
\begin{eqnarray*}
    |f(x-x_i)|\leq|f|(|x-x_i|)|&\leq& g(|x-x_i|)+(|f|-g)^{+}(|x-x_i|)
    \\&\leq&g(|x-x_i|)+(|f|-g)^{+}(|x|)+(|f|-g)^{+}(|x_i|)\\&\leq& \frac{\varepsilon}{3}+\frac{\varepsilon}{3}+\frac{\varepsilon}{3}=\varepsilon.
\end{eqnarray*}
Hence, for each $n\in \mathbb{N}$ and each $x\in A$ we have

   $$ |f_{n}(x)|\leq|f_{n}(x-x_i)|+|f_{n}(x_i)|\leq\varepsilon +\sup_{1\leq i\leq m}|f_{n}(x_i)|$$It follows that
    $\sup_{x\in A}|f_{n}(x)|\leq\varepsilon +\sup_{1\leq i\leq m}|f_{n}(x_i)|$. Since $f_{n}\xrightarrow{w^*}0$, we have $\sup_{1\leq i\leq m}|f_{n}(x_i)|\xrightarrow{n\to{\infty}}0$. Therefore, $\sup_{x\in A}|f_{n}(x)|\to 0$. That is, $(f_n)$ converges to zero uniformly on $A$. Proposition \ref{d} implies that $A$ is super weakly precompact.
\end{proof}
\begin{remark}
(1) In view of Theorem \ref{wo de PL-compact}, we have an alternative characterization of $|\sigma|(E,E^\prime)$--totally bounded sets in $E$. That is,  a subset $A$ of a Banach lattice $E$ is $|\sigma|(E,E^\prime)$--totally bounded if and only if  $T(A)\subset c_0$ is relatively compact for every almost Dunford-Pettis operator $T:E\to{c_0}$, or equivalently, if and only if  every weak$^{\,*}$-null almost $L$\,-sequence $(f_n)$ in $E^{\,\prime}$ converges uniformly to zero on $A$. This characterization helps know more about $|\sigma|(E,E^\prime)$--totally bounded sets and related operators. For instance, Theorem \ref{k.115}, Corollary \ref{u2y.3}, Theorem \ref{2y7.1}  improve the results  in Example 4.8 of \cite{Dodds and Fremlin PL-compact}. Nevertheless, in the sequel we use the term \textit{$|\sigma|(E,E^\prime)$--totally bounded} rather than \textit{super weakly precompact} in accord with the notion in the literature.
\par (2) From Corollary \ref{u2y.3} and Proposition \ref{image under order bounded operator} it follows that a set $A\subset E$ is $|\sigma|(E,E^\prime)$--totally bounded if and only if, for each Banach lattice $F$ with the positive Schur property and each order bounded linear operator $Q:E\to F$, the image $Q(A)$ is relatively compact.
\end{remark}
\par Now we turn our attention to the solid hull of a $|\sigma|(E,E^\prime)$--totally bounded set.  We can see that the solid hull of a $|\sigma|(E,E^\prime)$--totally bounded set is not necessarily  $|\sigma|(E,E^\prime)$--totally bounded. For instance, $B_{C[0,1]}=[-\mathbf{1},\mathbf{1}]=Sol\lbrace \mathbf{1}\rbrace$ is not $|\sigma|(C[0,1],C[0,1]^\prime)$--totally bounded.  Recall that an  operator $T:E\to{X}$ is called an \textit{AM-compact} (resp. \textit{order weakly compact}) \textit{operator} if  $T[-x,x]$ is relatively compact (resp. weakly compact) in $X$ for all $x\in E^+$.

\begin{theorem}\label{oh2y.3} For a Banach lattice $E$ the following statements are equivalent.
	\begin{itemize}
		\item [(a)] There holds $|f_n|\xrightarrow{w^*}0$ for every  almost $L$-sequence $(f_n)$ in $E^{\,\prime}$ such that $f_n\xrightarrow{w^*}0$.
		\item [(b)] The solid hull of every weakly precompact set in $E$ is $|\sigma|(E,E^\prime)$--totally bounded.
        \item [(c)] The solid hull of every $|\sigma|(E,E^\prime)$--totally bounded set in $E$ is also $|\sigma|(E,E^\prime)$--totally bounded.
		\item [(d)]	Every order interval of $E$ is $|\sigma|(E,E^\prime)$--totally bounded, that is,  every almost Dunford-Pettis operator $T:E\rightarrow c_0$ is \text{AM}-compact.
	\end{itemize}
\end{theorem}

\begin{proof} Only the implication $(a)\Rightarrow (b)$ needs a proof. Assume that $A\subset E$ is a weakly precompact set such that the solid hull $Sol(A)$ is not $|\sigma|(E,E^\prime)$--totally bounded. Then, by Proposition \ref{d}, there exist  a weak$^*$-null almost $L$-sequence $(f_n)$ in $E^\prime$ and a sequence $(z_n)\subset Sol(A)$ such that $\varepsilon<|f_n(z_n)|\le |f_n|(|z_n|)$ for some $\varepsilon>0$ and for all $n$. Therefore, for each $n$, there exist $x_n\in A$ such that $|z_n|\le |x_n|$ and hence, $\varepsilon< |f_n|(|z_n|)<|f_n|(|x_n|) $. By the Riesz-Kantorovich formula, for
	each $n$, there exists $g_n \in E^\prime$ such that $|g_n| \le | f_n|$ and $|g_n(x_n)| > \varepsilon$. By hypothesis (a), we have  $|f_n|\xrightarrow{w^*}0$ and so $|g_n|\xrightarrow{w^*}0$. Hence, from \cite[Lemma 2.2]{XCL}  it follows that the sequence $(g_n)$ is an $L$-sequence. However,
	 the weak precompactness of $A$ implies that $\varepsilon<|g_n(x_n)|\le \sup_{x\in A}|g_{n}(x)|\rightarrow0$ \cite[Corollary 9]{Ioana}. This is absurd.
\end{proof}
It should be noted that every  positively limited (resp. disjointly weakly compact) subset of a Banach lattice $E$ is weakly precompact if and only if every order interval in $E$ is weakly precompact. See \cite[Corollary 3.3]{positively limited set} and \cite[Theorem 2.4]{XCL} for details. The similar thing is true for  $|\sigma|(E,E^\prime)$--totally bounded sets.

\begin{corollary}\label{kgy.3} For a Banach lattice $E$ the following assertions are equivalent.
	\begin{itemize}
		\item [(a)] Each disjointly weakly compact set in $E$ is $|\sigma|(E,E^\prime)$--totally bounded.
		\item [(b)] Each almost Dunford-Pettis set in $E$ is $|\sigma|(E,E^\prime)$--totally bounded.
		\item [(c)]	Each positively limited set in $E$ is $|\sigma|(E,E^\prime)$--totally bounded.
		\item [(d)] Each order interval of $E$ is $|\sigma|(E,E^\prime)$--totally bounded.
	\end{itemize}
\end{corollary}

\begin{proof} $(a)\Rightarrow (b)\Rightarrow (c)$. This follows  from the facts that every positively limited set is an almost Dunford-Pettis set, and the latter is disjointly weakly compact. See \cite[Theorem 2.8 \& 3.2]{positively limited set} and  \cite[Remark 2.4(1)]{disjointly}.

	\par $(c)\Rightarrow (d)$. Follows easily from the observation that every order interval in a Banach lattice is  positively limited.
	
    \par $(d)\Rightarrow (a)$. Let $A\subset E$ be a disjointly weakly compact set. Since, by our hypothesis, every order interval of $E$ is $|\sigma|(E,E^\prime)$--totally bounded, and hence is certainly weakly precompact, it follows from \cite[Theorem 2.4]{XCL} that $A$ is weakly precompact. Hence, by Theorem \ref{oh2y.3}, $A$ is $|\sigma|(E,E^\prime)$--totally bounded.
\end{proof}

\begin{remark}\label{order interval is absolutely weakly totally bounded }
We can say a little more about a Banach lattice $E$ for which every order interval is  $|\sigma|(E,E^\prime)$--totally bounded. Recall that the absolute weak$^*$ topology $|\sigma|(E^{\prime},E)$ on $E^\prime$ is the locally convex-solid topology generated by the family of lattice seminorms $\{p_{x}:x\in E\}$, where $p_{x}(f)=|f|(|x|), \, f\in E^\prime$. It is should be noted that the absolute weak$^*$ topology $|\sigma|(E^{\prime},E)$ on $E^\prime$ is a Hausdorff, Lebesgue (i.e., order continuous), Levi topology and $E^\prime$ is $|\sigma|(E^{\prime},E)$--complete \cite[Theorem 6.5]{Locally Solid}.

\par (1) From a result of A. Grothendieck it follows that every order interval in $E$ is $|\sigma|(E,E^\prime)$--totally bounded if and only if every order interval in $E^\prime$ is $|\sigma|(E^{\prime},E)$--totally bounded. See, e.g., \cite[Theorem 3.27, Theorem 3.55; Exercise 8, p.180 ]{Positive}
\par(2) It is well known that a subset $A$ of a Hausdorff topological vector space is compact if and only if $A$ is totally bounded and complete.  Note that every order interval in $E^\prime$ is $|\sigma|(E^{\prime},E)$-closed and hence is $|\sigma|(E^{\prime},E)$-complete since $E^\prime$ is $|\sigma|(E^{\prime},E)$--complete. Therefore, it follows that every order interval in $E$ is $|\sigma|(E,E^\prime)$--totally bounded if and only if every order interval in $E^\prime$ is $|\sigma|(E^{\prime},E)$--compact, and this, in turn, is equivalent to that $E^\prime$ is discrete \cite[Corollary 6.57]{Locally Solid}.
\par(3) If every order interval in $E$ is $|\sigma|(E,E^\prime)$--totally bounded, then by Corollary \ref{kgy.3} every disjointly weakly compact set in $E$ is $|\sigma|(E,E^\prime)$--totally bounded and hence, from Theorem \ref{2y7.1} it follows that $E$ has weakly sequentially continuous lattice operations. That is, if $E^\prime$ is discrete, then $E$ has weakly sequentially continuous lattice operations.

\end{remark}

\par We are now in a position to characterize the weak Dunford-Pettis property  and the positive DP$^*$ property of Banach lattices in terms of $|\sigma|(E,E^\prime)$--totally bounded sets. The positive DP$^*$ property was recently introduced by the authors in \cite{positively limited set}. It was proved that a Banach lattice $E$ has the positive DP$^*$ property if and only if each disjointly weakly compact set in $E$ is positively limited \cite[Theorem 3.10]{positively limited set}.
\begin{theorem}\label{gy.3} For a Banach lattice $E$ the following statements are equivalent.
	\begin{itemize}
		\item [(a)] Every weakly precompact set in $E$ is  almost Dunford-Pettis  (resp. positively limited).
		\item [(b)] Every $|\sigma|(E,E^\prime)$--totally bounded set in $E$ is almost Dunford-Pettis  (resp. positively limited).
		\item [(c)]	$E$ has  the weak Dunford-Pettis property  (resp. positive DP$^{\,*}$ property).
	\end{itemize}
\end{theorem}

\begin{proof} $(a)\Rightarrow (b)$. Obvious.
	\par $(b)\Rightarrow (c)$. First,  for the weak Dunford-Pettis property, by \cite[Theorem 2.7]{Almost DP set} we have to show that $f_n(x_n)\rightarrow 0$, for every disjoint weakly null sequence $(x_n) \subset E$ and every disjoint weakly null sequence $(f_n)$ in $E^\prime$. Since $(x_n)$ is $|\sigma|(E,E^\prime)$--totally bounded, by hypothesis (b) $(x_n)$ is almost Dunford-Pettis, and so $f_n(x_n)\rightarrow 0$. Hence $E$ has the weak Dunford-Pettis property.
 \par As for the proof of the positive DP$^{\,*}$ property,  in view of Theorem 3.10 of \cite{positively limited set} it suffices to show that $f_n(x_n)\rightarrow0$ holds for every disjoint weakly null sequence $(x_n) \subset E_{+}$ and every weak$^*$-null sequence $(f_n)$ in $(E^\prime)^{+}$. Again, by hypothesis (b), the $|\sigma|(E,E^\prime)$--totally boundedness of the sequence $(x_n)$ implies that $(x_n)$ is positively limited, and hence $f_n(x_n)\rightarrow 0$.

\par	$(c)\Rightarrow (a)$. This follows directly from  \cite[Theorem 2.9]{disjointly} (resp. \cite[Theorem 3.10]{positively limited set}) since every weakly precompact set in a Banach lattice is necessarily disjointly weakly compact (cf. \cite[Proposition 2.5.12 iii)]{Meyer}).
\end{proof}

\section{PL-compact operators in Banach lattices}

Let $E$ and $F$ be  Banach lattices. If $x\in E^+$, we denote by $E_x$ the order ideal generated by $x$ in $E$ (equipped with the sup norm) and denote by $i_x$ the injection of $E_x$ into $E$. If $g\in (F^{\,\prime})^+$, $j_g$ and $(F;g)$ are as in the preceding section. P.G Dodds and D.H. Fremlin introduced the notions of AMAL-compact operators and PL-compact operators to  study  compact operators. For our convenience we list the  definitions.

\begin{definition}\cite{Dodds and Fremlin PL-compact}\label{2h.1}
   (1)  An operator  $T\in \mathcal{L}(X,E)$ is called PL-compact if the composition $j_{g}T:X\to{(E;g)}$ is a compact operator for each $g\in (E^{\,\prime})^+$, i.e.,  if $T(B_X)$ is a $|\sigma|(E,E^\prime)$--totally bounded set in $E$.

   \par (2) An operator  $T\in \mathcal{L}(E,F)$  is  called AMAL-compact if the bicomposition $j_{g}Ti_{x}:E_x\to(F;g)$ is compact for each $x\in E^+$ and each $g\in (F^{\,\prime})^+$, i.e., if $T[-x,x]$ is  $|\sigma|(F,F^\prime)$--totally bounded in $F$ for each $x\in E^+$.
 \end{definition}
\begin{remark}\label{wo de remark3}
 (1) In view of Proposition \ref{d} and Theorem \ref{wo de PL-compact}, a bounded linear operator $T:X\rightarrow E$ is PL-compact  if and only if we have $\|T^{\prime}f_{n}\|\to{0}$ for each weak$^*$-null almost $L$-sequence $(f_{n})\subset E^{\prime}$. Similarly, an operator $T:E\to{F}$  is  AMAL-compact  if and only if  $|T^{\prime}f_{n}|\stackrel{w^{*}}{\longrightarrow}0$ for each weak$^*$-null almost $L$-sequence  $(f_{n})\subset F^{\prime}$.

 \par (2) The name \textit{PL-compact operator} is justified by the observation that  $T:X\to E$ is PL-compact  if and only if for each $L$-space $F$ and each positive operator $Q:E \to F$ the composition $QT$ is a compact operator \cite[3.7.E1]{Meyer}. Also, Theorem 4.9 of \cite{Dodds and Fremlin PL-compact} asserts that an operator  $T\in \mathcal{L}(X,E)$ is   PL-compact  if and only if $T^{\,\prime}$ is AM-compact.
 \par (3) An operator  $T\in \mathcal{L}(E,F)$ is   AMAL-compact  if and only if $T^{\,\prime}[-g,g]$ is $|\sigma|(E^{\prime},E)$--totally bounded for each $g\in (F^{\prime})^+$ (cf. e.g., \cite[Exercise 8, p.180]{Positive}). Also, a positive operator dominated by a positive AMAL-compact operator is likewise AMAL-compact \cite[Theorem 5.11]{Positive}.
\end{remark}

We can see that  an operator  $T\in \mathcal{L}(X,E)$ maps relatively weakly compact subsets of $X$ onto $|\sigma|(E,E^{\,\prime})$--totally bounded sets in $E$ if and only if the composite $j_{f}T:X\to (E;f)$ is a Dunford-Pettis operator for each $f\in (E^{\prime})^+$, that is, $|Tx_n|\xrightarrow{w}0$ holds in $E$  whenever $x_{n}\xrightarrow{w}0$ in $X$.
 The next result tells us when an operator maps  disjointly weakly compact sets onto totally bounded sets.
\begin{theorem}\label{send dwc to totally bounded}
An operator $T\in \mathcal{L}(E,F)$ sends each disjointly weakly compact set to a $|\sigma|(F,F^{\,\prime})$--totally bounded set if and only if $T$ is AMAL-compact and  $|Tx_n|\xrightarrow{w}0$ in $F$ for  every (disjoint) sequence $(x_{n})\subset E$ satisfying $x_{n}\xrightarrow{w}0$.
\par For an order bounded operator $T:E\to F$, if $T$ is AMAL-compact, then $T$ maps disjointly weakly compact sets  onto $|\sigma|(F,F^{\,\prime})$--totally bounded sets.
\end{theorem}

\begin{proof}
\par  Assume that $T$ carries each disjointly weakly compact subset of $E$ to a $|\sigma|(F,F^{\,\prime})$--totally bounded subset of $F$. Then, for each $g\in (F^{\prime})^+$, $j_{g}T:E\to (F;g)$ maps each disjointly weakly compact subset of $E$ onto a relatively compact subset of $(F;g)$. Thus, $j_{g}T:E\to (F;g)$ is both AM-compact and Dunford-Pettis. That is, $T$ is AMAL-compact, and for  every sequence $(x_{n})\subset E$ satisfying $x_{n}\xrightarrow{w}0$, we have $g|Tx_n|=\|j_{g}T(x_n)\|\to 0$.
\par For the converse, assume that $T:E\to F$ be an AMAL-compact operator such that $|Tx_n|\xrightarrow{w}0$ in $F$ for  every disjoint weakly null sequence $(x_{n})$ of $E$. That is,  for each $g\in (F^{\prime})^+$,  $j_{g}T:E\to (F;g)$ is both AM-compact and almost Dunford-Pettis.
Now, let $A$ be a disjointly weakly compact subset of $E$ and let $g\in (F^{\prime})^+$ be fixed. By definition, we can assume without loss of generality that $A$ is solid. Then, for every disjoint sequence $(x_n)\subset A$, we have $\|j_{g}Tx_n\|\to 0$ since $j_{g}T$ is almost Dunford-Pettis. Then in view of Theorem 4.36 of \cite{Positive}, for each $\varepsilon>0$  there exists some $u\in{E^{+}}$  such that
	$\|j_{g}T[(|x|-u)^{+}]\|<\frac{\varepsilon}{2}$ holds for all $x\in A$. Therefore, from the identity $|x|=|x|\wedge u+(|x|-u)^+$ we can see that
$$j_{g}T(A)\subset j_{g}T[-u,u]+\varepsilon B_{(F;g)}$$
Since $j_{g}T$ is AM-compact,  it follows that  $j_{g}T(A)$ is relatively compact. This implies that $T(A)$ is $|\sigma|(F,F^{\,\prime})$--totally bounded.

\par For an order bounded operator $T\in \mathcal{L}^{b}(E,F)$, the AMAL-compactness of $T$ implies that $|Tx_n|\xrightarrow{w}0$ holds in $F$  whenever $x_{n}\xrightarrow{w}0$ in $E$. This is due to C.D. Aliprantis and O. Burkinshaw (see \cite{Aliprantis} or \cite[Theorem 5.96]{Positive}). Thus, the proof is finished by the first part.
\end{proof}

\par The identity operator $I:\ell^1\to\ell^1$ is an example of an AMAL-compact operator which is not PL-compact. Theorem 4.4 of \cite{Dodds and Fremlin PL-compact} asserts that if $E$ and $F$ are Banach lattices such that $E^\prime$ has order continuous norm, then every order bounded AMAL-compact operator $T:E\to F$ is PL-compact. The following result shows that the order continuity of the norm of $E$ is indeed a necessary and sufficient condition.

 \begin{corollary}
	For a Banach lattice $E$ the following statements are equivalent.
\begin{enumerate}
  \item [(a)] $E^{\,\prime}$  has order continuous norm.
  \item [(b)] Every order bounded  AMAL-compact  operator from $E$ to an arbitrary Banach lattice $F$ is PL-compact.
  \item [(c)] Every  AMAL-compact positive operator from $E$ to an arbitrary Banach lattice $F$ is PL-compact.
\end{enumerate}
\end{corollary}

\begin{proof}
\par(a)$\Rightarrow$(b)  is included in Theorem 4.4 of \cite{Dodds and Fremlin PL-compact}. Here, as an immediate consequence of Theorem \ref{send dwc to totally bounded}, the implication is obvious since $B_E$ is disjointly weakly compact.
\par (b)$\Rightarrow$(c) Obvious.
\par(c)$\Rightarrow$(a) Assume that  $E^{\,\prime}$ fails to have order continuous norm. Then, $E$ contains a closed sublattice $U$ lattice isomorphic to $\ell^1$ (cf. \cite[Theorem 2.4.14]{Meyer}). Furthermore, $U$ is the range of a positive projection $P$ on $E$. See \cite[Proposition 2.3.11]{Meyer}.  We can easily verify that the positive projection $P:E\to{U}$ is an AMAL-compact operator. However, $P$ is not PL-compact.
\end{proof}
Let us recall that, by Corollary \ref{u2y.3}, a Banach lattice $E$ has the positive Schur property if and only if every $|\sigma|(E,E^{\,\prime})$--totally bounded subset of $E$ is relatively compact. Next, we consider those operators which carry $|\sigma|(E,E^{\,\prime})$--totally bounded subsets of $E$ onto relatively compact sets.

\begin{proposition}\label{h3.4} For a Banach lattice $E$ the following statements are equivalent:
	\begin{itemize}
		\item [(a)] $E$ has the positive Schur property.
		\item [(b)] Every bounded linear operator from $E$ into an arbitrary Banach space $X$ carries $|\sigma|(E,E^{\,\prime})$--totally bounded sets  onto relatively compact sets.
		\item [(c)] Every bounded linear operator from $E$ into  $\ell_{\infty}$ carries $|\sigma|(E,E^{\,\prime})$--totally bounded sets  onto relatively compact sets..
	\end{itemize}
\end{proposition}

\begin{proof}
	
	It suffices to prove that $(c)\Rightarrow (a)$. Assume to the contrary that $E$ does not have the positive Schur property. Then there exists a weakly null and disjoint  sequence $(x_n)$
	in $E^+$ such that $\|x_n\|=1$ for all $n$. Choose a normalized
	sequence $(f_n)$ in $E^\prime$ such that $f_{n}(x_n)=1$ for all $n$. We define the operator $T:
	E\rightarrow \ell_{\infty} $ by
	$$Tx=(f_{n}(x) )~~~~~~~ , ~~~~~~~~x\in E.$$Note that $(x_n)$ is a disjoint weakly null sequence  (and hence $|\sigma|(E,E^{\,\prime})$--totally bounded set). However, $T$  does not map $(x_n)$ onto a relatively compact set  since  $\|Tx_n\|
	\geq1$  for all $n\in \mathbb{N}$. This leads to a contradiction.	
\end{proof}

\begin{theorem}\label{cc}For two Banach lattices $E$ and $F$,  if every PL-compact operator $T:E\to F$ is Dunford-Pettis, then one of the following two conditions holds:
	\begin{enumerate}
		\item [(a)] $E$ has the weakly sequentially continuous lattice operations.
		\item [(b)] $F$ is a KB-space.
	\end{enumerate}	
\end{theorem}

\begin{proof}
	Suppose that neither (a) nor (b) holds. Then, as in Theorem \ref{2y7.1}, we would find a weakly null sequence $(x_n)\subset E$ and an almost $L$-sequence $(h_{n})$ in $E^\prime$ with $h_{n}\xrightarrow{w^*}0$  such that  $h_n(x_n)\ge \varepsilon$ for all $n$.  Consider the PL-compact operator
	$S:E\rightarrow c_0$ defined by
	$Sx=(h_n(x))~~~~~~~ , ~~~~~~~~x\in E$. Since $F$ is not a KB-space, there is a lattice embedding $i: c_0
	\to F$.
	It is clear that the operator
	$T:=i~ o~S: E \to F$ is PL-compact,  but $T$ is not Dunford-Pettis since  $\|Tx_n\|\nrightarrow 0$ with $x_{n}\xrightarrow{w}0$.
\end{proof}

\begin{theorem}\label{Dunford is PL-compact} For two Banach lattices $E$ and $F$,  if every Dunford-Pettis  operator $T:E\to F$ is PL-compact, then one of the following two conditions holds:
	\begin{enumerate}
		\item [(a)] $E^\prime$ has order continuous norm.
		\item [(b)] $F^\prime$ has order continuous norm.
	\end{enumerate}	
\end{theorem}

\begin{proof}
	Suppose that neither (a) nor (b) holds. Then there would exist an order bounded disjoint  sequence $(\phi_n)\subset (E^\prime)^+$  such that $\|\phi_n\|=1$ for all $n$.  Consider the operator
	$U:E\rightarrow \ell_1$ defined by
	$U(x)=(\phi_n(x))~~~~~~~ , ~~~~~~~~x\in E$. Since $F^\prime$ does not have order continuous norm, there is exist an order bounded disjoint sequence $(f_n)\subset (F^\prime)^+$  such that $\|f_n\|=1$ for all $n$.  Hence for each $n$ there exists $y_n \in F^+$ with $\|y_n\|=1$ and
	$f_n(y_n)\ge \frac{1}{2}$. Define a positive operator $V:\ell_1\rightarrow F$ defined by
	$V(\lambda)=\sum_{n=1}^{\infty}\lambda_ny_n~~~~~~~ , ~~~~~~~~\lambda= (\lambda_n)\in \ell_1$.
	Consider the operator
	$T:= V o U: E\rightarrow \ell_1\to F$ defined by
	$Tx=\sum_{n=1}^{\infty}\phi_n(x)y_n~~~~~~~ , ~~~~~~~~x\in E$. It is clear that $T$ is
	Dunford-Pettis. However, $T$ is not PL-compact. To see this, note that if $h \in F^\prime$ then $T^{\prime}h=\sum_{n=1}^{\infty}h(y_n)\phi_n$. In particular, for every $k$ we have
	$\|T^{\prime}(f_k)\|\ge \|f_k(y_k) \phi_k\|=f_k(y_k) \ge \frac{1}{2}$. As $(f_k)$ is a weak$^*$-null and almost $L$-sequence  in $F^\prime$ (since it is an order bounded disjoint sequence), $T$ is not PL-compact.
\end{proof}

\par We now turn our attention  to the domination problem of  PL-compact operators. First we give an example to show that a positive operator dominated by a PL-compact operator is not necessarily PL-compact.

\begin{example}\label{composition counterexample}
	Let $0\le S\le T:   \ell_{1}\rightarrow L^1[0,1]$ be defined as
	\begin{eqnarray*}
		S(\alpha) &=& \sum_{n=1}^{\infty}\alpha_nr_n^+~~~~~~~ , ~~~~~~~~\alpha=(\alpha_n)\in \ell_1\\
		T(\alpha) &=& (\sum_{n=1}^{\infty}\alpha_n)\cdot\textbf{1}~~~~~~~ , ~~~~~~~~\alpha=(\alpha_n)\in \ell_1
	\end{eqnarray*}
	where $r_{n}(t)$ is the $n$th Rademacher function on $[0, 1]$.
	$T$ is compact and hence PL-compact. However, $S$ is not PL-compact. %Indeed, $(Se_n)_{n=1}^{\infty}=(r_n^+)$ is not strong weak limited in $L^1[0,1]$.
\end{example}

 Theorem 4.6 of \cite{Dodds and Fremlin PL-compact} asserts that, for two Banach lattices $E$ and $F$ such that $E^{\,\prime}$ and $F$ have order continuous norms, the set of regular PL-compact operators from $E$ to $F$ forms a band in $\mathcal{L}^{r}(E,F)$. It should be noted that, in a Banach lattice $E$ with order continuous norm, the norm topology and $|\sigma|(E,E^{\,\prime})$ agree on every order interval of $E$, and hence have the same order bounded totally bounded sets. See, e.g., \cite[Theorem 4.17]{Positive}. Therefore, if $E$ has order continuous norm but $E$ is not discrete, then there exists an order interval of $E$ which is not $|\sigma|(E,E^{\,\prime})$-totally bounded. On the other hand, a Banach lattice for which every order interval is totally bounded with respect to the  absolute weak topology does not necessarily have order continuous norm. For instance, by Theorem \ref{k.115},  $B_{c}=[-\mathbf{1},\mathbf{1}]$ is $|\sigma|(c,c^{\,\prime})$--totally bounded. The next result tells us when every  positive operator dominated by a PL-compact operator is always PL-compact.

\begin{theorem}\label{domination} For two Banach lattices $E$ and $F$  the following assertions are equivalent:
	\begin{enumerate}
		\item Every positive operator from $E$ into $F$ dominated by a PL-compact operator is PL-compact.
		\item One of the following two conditions holds:
		\begin{enumerate}
			\item [(a)] The norm of $E^{\,\prime}$ is order continuous.
			\item [(b)] Every order interval in $F$ is $|\sigma|(F,F^{\,\prime})$--totally bounded.	
		\end{enumerate}	
	\end{enumerate}
\end{theorem}

\begin{proof} $(1)\Rightarrow(2)$ Assume by way of contradiction that  the norm of $E^\prime$ is not order continuous and  there exists an order interval $[-y,y]\subset F$ ($y\in F^+$) which is not $|\sigma|(F,F^{\,\prime})$--totally bounded. To finish the proof, we have to construct two positive operators $0\leq S \leq T:E\to{F}$ such that $T$ is PL-compact and $S$ is not PL-compact.

\par Since the norm of $E^{\prime}$ is not order continuous,  there  exists some $x^{\prime}\in (E^\prime)^+$ and a disjoint sequence $(x_{n}^\prime)\subset[0, x^{\prime}]$  such that $\|x_{n}^\prime\|=1$ for all $n$. Let us define a positive operator $S_{1}:E\to{\ell_{1}}$ by
$S_{1}x=(x_{n}^{\prime}(x))_{n=1}^{\infty}$
for all $x\in{E}$.  Also, since the order interval $[-y,y]\subset F$  is not $|\sigma|(F,F^{\,\prime})$--totally bounded, by Proposition \ref{d} we can find  a weak$^{*}$-null almost $L$\,-sequence $(y_{n}^\prime)\subset{F^{\,\prime}}$ and a sequence $(y_n)\subset [0,y]$ such that ${\vert{{y_{n}^\prime}(y_n)}\vert}>{\epsilon}$ for some $\epsilon>{0}$ and all ${n}\in{\mathbb{N}}$.  We define the positive operator $S_{2}:{\ell_{1}}\to{F}$  by ${S_{2}((\lambda_{n}))}={\sum_{n=1}^{\infty}{\lambda}_{n}{y_n}}$ for all ${({\lambda}_{n}})\in{\ell_{1}}$.

\par Now consider two  positive operators $S,T:E\to{F}$, defined by
	\begin{center}
		$S(x)=S_{2}S_{1}(x)={{\sum_{n=1}^{\infty}}{x_{n}^{\prime}(x)}{y_n}}$\,\, and \,\,\,$T(x)=x^{\prime}(x)y$.
	\end{center}
for all $x\in E$. It is easy to verify that $0\le{S}\le{T}$. Clearly, $T$ is compact, and hence $T$ is PL-compact. However,
	\begin{center}
		$\|S^{\prime}y_{n}^{\prime}\|\ge{{\vert{y_{n}^{\prime}(y_n)}\vert}{\|x_{n}^\prime}\|}>{\epsilon}$
	\end{center}
	for all $n$. This implies that $S$ is not PL-compact (see Remark \ref{wo de remark3} (1)).

\par $(2)\, (a)\Rightarrow(1)$ Let $S, T:E\rightarrow F$ be two positive operators  with $0\leq S\leq T$ and $T$  PL-compact. Then we have $0\leq S^{\,\prime}\leq T^{\,\prime}:F^{\,\prime}\rightarrow E^{\,\prime}$. Since $T$ is PL-compact,  $T^{\,\prime}$ is AM-compact \cite[Theorem 4.9 ]{Dodds and Fremlin PL-compact}. If the norm of $E^{\,\prime}$ is order continuous, then by a result of Dodds and Fremlin \cite[Theorem 4.7]{Dodds and Fremlin PL-compact} we know that $S^{\,\prime}$ is also AM-compact, that is, $S$ is PL-compact.

\par $(2)\,(b)\Rightarrow(1)$ Let $S, T:E\rightarrow F$ be two positive operators  such that $0\le S\le T$ and $T$ is  PL-compact. If every order interval in $F$ is $|\sigma|(F,F^{\,\prime})$--totally bounded, then by Theorem \ref{oh2y.3} $Sol(TB_E)$ is an $|\sigma|(F,F^{\,\prime})$--totally bounded subset of $F$. Hence, from the set inclusin $SB_E\subset Sol(TB_E)$ we can easily see that $SB_E$ is $|\sigma|(F,F^{\,\prime})$--totally bounded; that is, $S$ is also  PL-compact.
\end{proof}

\bibliographystyle{amsplain}

\end{document}